\documentclass[11pt]{amsart}
\usepackage{amssymb,amsmath,amsthm,amsfonts}
\usepackage{mathrsfs,a4wide,dsfont}
\usepackage[T1]{fontenc}
\usepackage{color}
\usepackage{ifpdf}
\ifpdf
    \usepackage[pdftex]{graphicx}
    \DeclareGraphicsExtensions{.png,.pdf,.jpg}
\else
   \usepackage[dvips]{graphicx}
   \DeclareGraphicsExtensions{.eps}
\fi

\theoremstyle{plain}
\newtheorem{theorem}{Theorem}[section]

\newtheorem{corollary}[theorem]{Corollary}
\newtheorem{remark}[theorem]{Remark}

\theoremstyle{definition}
\theoremstyle{remark}
\numberwithin{equation}{section}
\newenvironment{oss}{\begin{remark} \begin{rm}}{\end{rm} \end{remark}}

\newcommand{\hs}{{\mathscr H}}

\newcommand{\R}{{\mathbb R}}

\newcommand{\N}{{\mathbb N}}
\newcommand{\Z}{{\mathbb Z}}

\newcommand{\Om}{\Omega}

\newcommand{\e}{\varepsilon}
\newcommand{\SM}{\mathcal{SM}}
\newcommand{\SMeC}{\SM_{\e, C}}
\newcommand{\esssup}{\mathop{\mathrm{ess\,sup}}}
\newcommand{\essinf}{\mathop{\mathrm{ess\,inf}}}
\newcommand{\loc}{\mathrm{loc}}
\renewcommand{\P}{\mathrm{Per}}
\newcommand{\PhC}{\P_{h_C}}
\newcommand{\dist}{\mathrm{dist}}
\newcommand{\distC}{\mathrm{dist_C}}

\newcommand{\wtos}{\stackrel{*}{\rightharpoonup}}
\newcommand{\ov}[1]{\overline{#1}}

\newcommand{\mrestr}{\ \rule{.4pt}{7pt}\rule{6pt}{.4pt}\ }

\let\TeXchi\chi
\newbox\chibox
\setbox0 \hbox{\mathsurround0pt $\TeXchi$}
\setbox\chibox \hbox{\raise\dp0 \box 0 }
\def\chi{\copy\chibox}

\title[The Anisotropic Outer Minkowski Content]{A Remark On The Anisotropic Outer Minkowski Content}

\author[A. Chambolle]
{Antonin Chambolle}
\address[A.\,Chambolle]{CMAP, \'{E}cole Polytechnique, CNRS, F-91128, Palaiseau, France}
\email[A. Chambolle]{antonin.chambolle@cmap.polytechnique.fr}
\author[S. Lisini]
{Stefano Lisini}
\address[S.\,Lisini]{Dipartimento di Matematica ``F.\,Casorati'', Universit\`a degli Studi di Pavia, via Ferrata 1, I-27100 Pavia, Italy}
\email[S. Lisini]{stefano.lisini@unipv.it}
\urladdr{http://www-dimat.unipv.it/lisini/}
\author[L. Lussardi]
{Luca Lussardi}
\address[L.\,Lussardi]{Dipartimento di Matematica e Fisica ``N.\,Tartaglia'', Universit\`a Cattolica del Sacro Cuore, via dei Musei 41, I-25121 Brescia, Italy}
\email[L. Lussardi]{l.lussardi@dmf.unicatt.it}
\urladdr{http://www.dmf.unicatt.it/~lussardi/}

\begin{document}
\begin{abstract}
\small{We study an anisotropic version of the outer Minkowski content of a closed set in $\R^n$.
In particular, we show that it exists on the same class of sets for which the classical outer Minkowski content coincides with the Hausdorff measure, and we give its explicit form.
\vskip .3truecm
\noindent Keywords: outer Minkowski content, finite perimeter sets
\vskip.1truecm
\noindent 2010 Mathematics Subject Classification: 28A75, 49Q15, 52A39}
\end{abstract}
\maketitle

\section{Introduction}
As it is well known, the classical {\it Minkowski content} of a closed set $S\subset \R^n$ is defined by
\begin{equation}\label{classMink}
\mathcal M(S)\ := \ \lim_{\e \to0^+}\frac{|\{x \in \R^n : \dist(x,S)\leq \e\}|}{2\e}
\end{equation}
whenever the limit in \eqref{classMink} exists and is finite; here $|\cdot|$ denotes the Lebesgue measure in $\R^n$. The quantity $\mathcal M$ measures the area of ``$(n-1)$-dimensional sets'', and it is an alternative to the more classical Hausdorff measure $\mathscr H^{n-1}$. With the role of surface measure, the Minkowski content turns out to be important in many problems arising from real applications: for instance $\mathcal M$ is related to evolution problems for closed sets \cite{A,L,SZ}.
\par
Clearly, it poses as natural problems its existence and comparison with $\hs^{n-1}$. Let us mention some known results in this direction. In \cite[p.~275]{Federer} the author proves that $\mathcal M(S)$ exists and coincides
with $\hs^{n-1}(S)$ whenever $S$ is compact and {\it $(n-1)$-rectifiable}, i.e.\,$S=f(K)$ for some $K\subset \R^{n-1}$ compact and $f\colon \R^{n-1}\to\R^n$ is Lipschitz.
A generalization of this result is contained in \cite[p.~110]{AFP-libro}. Here,  the authors consider {\it countably $\hs^{n-1}$-rectifiable} compact sets in $\R^n$, i.e.\,sets which can be covered by a countable family of sets $S_j$, with $j \in \N$, such that $S_0$ is $\hs^{n-1}$-negligible and $S_j$ is a $(n-1)$-dimensional surface in $\R^n$ of class $C^1$, for any $j>0$.
In this case, $\mathcal M(S)$ exists and coincides with $\hs^{n-1}(S)$ if a further density assumption on $S$ holds: more precisely there must exist $\gamma>0$ and $\eta$ a probability measure on $\R^n$ satisfying $\eta(B(x,r))\geq \gamma r^n$, for each $r\in (0,1)$ and for each $x\in S$, where $B(x,r)$ is the open ball centered in $x$ of radius $r$.
Counterexamples \cite[p.~109]{AFP-libro} show that the countable rectifiability is indeed not sufficient to ensure the existence of $\mathcal M$.
\par
More recently in \cite{ACV}, motivated by problems in stochastic geometry, a generalization of the Minkowski content has been introduced, the so-called {\it outer Minkowski content}, which is defined by
\begin{equation}\label{classSMink}
\SM(E)\ := \ \lim_{\e \to0^+}\frac{|\{x \in \R^n : \dist(x,E)\leq \e\}\setminus E|}{\e}, \quad E \subset \R^n\,\,{\rm compact}.
\end{equation}
In \cite{ACV} the authors investigate general conditions ensuring the existence of $\SM$: in particular, they prove that if $E$ is a set with finite perimeter and $\mathcal M(\partial E)$ exists and coincides with the perimeter of $E$, then also $\SM(E)$ exists and coincides with the perimeter of $E$ (in $\Om$).
\par
Now, notice that the quantity which appears in the argument of the limit in \eqref{classSMink} can be rewritten as (provided the set $E$ is ``nice'' enough)
$$
\frac{1}{\e}(|E+\e B(0,1)|-|E|).
$$

We consider in this short note a variant of this content, which is
an anisotropic outer Minkowski content. The idea is to study the limit,
as $\e \to 0^+$, of
\begin{equation}\label{anisq}
\frac{1}{\e}(|E+\e C|-|E|),
\end{equation}
where $C\subset \R^n$ is a closed convex body.
It is standard that if $E$ is convex, then
$|E+\e C|$ is a polynomial in $\e$ (of degree $n$) whose coefficient of the first degree term (see also Remark \ref{remConv} below)
is precisely the anisotropic perimeter
\begin{equation}\label{anisoper}
\int_{\partial E}h_C(\nu_E)\,d\mathscr H^{n-1}\,,
\end{equation}
where $h_C$ is
the support function of $C$, defined by $h_C(\nu)=\sup_{x\in C}x\cdot\nu$, and $\nu_E$ the outer normal to $\partial E$,
see~\cite{Schneider} for details. The convergence of~\eqref{anisq}
to~\eqref{anisoper} follows for convex sets $E$ and can
be easily extended to (very) smooth sets.

We show here
two (expected) results: first, as a functional defined on sets,
\eqref{anisq} $\Gamma$-converges to the natural limit~\eqref{anisoper} as $\e\to 0$.

Second, we show in Theorem \ref{ThPW} that given any set
for which the
(classical) outer Minkowski content equals the perimeter, then the limit
of \eqref{anisq}, as $\e\to 0^+$, coincides with~\eqref{anisoper}.

The proof of Theorem \ref{ThPW} is quite technical, because we wanted
to work under the only assumption of the convergence of the classical
content. We show that this convergence implies that the boundary
is flat enough in a relatively uniform way, so that the convergence
of \eqref{anisq} holds.
\par

Eventually, we also deduce a $\Gamma$-convergence result (see \cite{B,DM-libro} for details) for  functionals of the type
\begin{equation}\label{anisq1}
\frac{1}{\e}\int\big( \esssup_{x-\e C} u-u(x)\big)\,dx
\end{equation}
which coincides with \eqref{anisq} on characteristic functions of sets.
The limit is (quite obviously) given by $\int_\Om h_C(-Du)$ (where
the minus sign accounts for the fact that the outer normal was appearing
in \eqref{anisoper}, and not the inner normal which corresponds more naturally
to the gradient of the characteristic function $\chi_E$)
\smallskip

As a simple
corollary, one is able to approximate functionals of the type
$$
\int_{\partial E}\phi(\nu_E)\,d\hs^{n-1},
$$
for $\phi$ a positively one-homogenous convex function $\phi:\R^n\to [0,+\infty)$ (and positive away from $0$).
Indeed, it suffices to choose the convex body
$$
C\ :=\ \{x \in \R^n : x \cdot \nu \le \phi(\nu) \ \forall\, \nu \in \R^n\}
$$
and apply our results.
\\
\par
The paper is organized as follows: in section \ref{secStatements} we define the setting and we state the results, then in section \ref{secProofs} we prove the $\Gamma$-convergence result for \eqref{anisq1}, and then the pointwise convergence result for \eqref{anisq}.

\section{Notation and preliminaries}

\subsection{Notation}
Let $n\geq 1$ be integer. Given a measurable set $A \subset \R^n$, we will denote by $|A|$ its Lebesgue measure. If $k\in \{0,\dots,n\}$, the $k$-dimensional Hausdorff measure of $S \subset \R^n$ will be denoted by $\hs^k(S)$. We will use the notation $x\cdot y$ for the standard scalar product in $\R^n$ between $x$ and $y$, $B(x,r)$ for the closed ball of radius $r$ centered in $x$. Finally, here convergence in $L^1_\loc(\Om)$ means convergence in $L^1(B\cap \Om)$ for any ball $B$. (Strictly speaking, it is thus the convergence in $L^1_\loc(\R^n)$
of the functions extended with the value $0$ outside of $\Om$.)

We say that a sequence of sets $E_j\subset \R^n$ converges to $E\subset \R^n$ in $L^1_\loc(\Om)$ as $j\to +\infty$,
if $\chi_{E_j}$ converges to $\chi_E$ in $L^1_\loc(\Om)$ as $j\to +\infty$, where $\chi_S$ denotes the characteristic function of the set $S$.

\subsection{Geometric measure theory}

In this paragraph we recall some basic notions of geometric measure theory we will need; for an exhaustive treatment of the subject we refer the reader to \cite{S}.

Let $n\geq 1$ be integer and let $k \in \N$ with $k \leq n$. We say that $S\subset \R^n$ is {\it countably $\hs^k$-rectifiable} if $S$ can be covered by a countable family of sets $S_j$, with $j \in \N$, such that $S_0$ is $\hs^k$-negligible and $S_j$ is a $k$-dimensional surface in $\R^n$ of class $C^1$, for any $j>0$.

Let $E\subset \R^n$ be a measurable set and $\Omega\subset \R^n$ be an open domain. We say that $E$ has {\it finite perimeter in $\Omega$} if the distributional derivative of $\chi_E$, denoted by $D\chi_E$, is a $\R^n$-valued Radon measure on $\Om$ with finite total variation; the perimeter of $E$ in $\Omega$ is defined by $\P(E;\Om):=|D\chi_E|(\Om)$, where $|D\chi_E|$ denotes the total variation of $D\chi_E$.
The {\it upper} and {\it lower $n$-dimensional densities of $E$ at $x$} are respectively defined by
$$
\Theta_n^*(E,x):=\limsup_{r \searrow 0}\frac{|E\cap B(x,r)|}{\omega_n r^n}, \quad \Theta_{*n}(E,x):=\liminf_{r \searrow 0}\frac{|E\cap B(x,r)|}{\omega_n r^n},
$$
where $\omega_n$ is the volume of the $n$-dimensional unit ball. If $\Theta_n^*(E,x)=\Theta_{*n}(E,x)$ their common value is denoted by $\Theta_n(E,x)$.
For every $t\in [0,1]$ we define
\begin{equation}\label{defEt}
E^t:=\{x\in\R^n : \Theta_n(E,x)=t\}.
\end{equation}
The {\it essential boundary of $E$} (or measure-theoretic boundary)
is the set $\partial_* E:=\R^n \setminus (E^0\cup E^1)$. It turns out that if $E$ has finite perimeter in $\Om$, then
$\hs^{n-1}(\partial_* E\setminus E^{1/2})=0$, and $\P(E;\Om)=\hs^{n-1}(\partial_* E \cap \Om)$.

Moreover, one can define a subset
of $E^{1/2}$ as the set of points $x$
where there exists a unit vector $\nu_E(x)$ such that:
\[
\frac{E-x}{r} \to \{y \in \R^n: y\cdot \nu_E(x) \le 0\},
\quad\text{ in $L^1_{\rm loc}(\R^n)$ as $r \to 0$},
\]
and which is referred to as the {\it outer normal to $E$ at $x$}.
The set where $\nu_E(x)$ exists is called the \textit{reduced boundary}
and is denoted by  $\partial^* E$.
One can show that $\hs^{n-1}(\partial_* E\setminus \partial^* E)=0$, moreover,
one has the decomposition
$D\chi_E\ =\ (-\nu_E)\hs^{n-1}\mrestr\partial^* E$.

\section{Statement of the results}\label{secStatements}

Let us assume that $C\subset\R^n$ is a closed
convex body, that is, bounded and with $0$ in its
interior. We denote its support function by
$h_C(\nu)\, =\, \sup_{x\in C}x\cdot \nu$,
and its polar function  is $h^\circ_C(x)\,:=\,\sup_{h_C(\nu)\le 1} x\cdot \nu$.
It is well known, then, that both $h_C$ and $h^\circ_C$ are convex, one-homogeneous and Lipschitz functions, moreover $C=\{ h^\circ_C\le 1\}$.

By assumptions, there also exist $a,b$ with $0<a<b$ such that
$B(0,a)\subseteq C\subseteq B(0,b)$, in particular, we have for
all $\nu,x\in \R^n$
\begin{equation}\label{comparh}
a|\nu|\, \le\, h_C(\nu)\,\le\, b|\nu|\,,\qquad \frac{1}{b}|x|\,\le\,
h^\circ_C(x)\,\le\, \frac{1}{a}|x|\,.
\end{equation}

Let $\Om\subset \R^n$ be an open domain.
Given a Lebesgue measurable set $E\subset \Om$, we introduce the outer $\e,C$-Minkowski content,
\begin{equation}\label{defsmec0}
\SMeC^0(E;\Om)\ :=\ \frac{1}{\e} (|\Om\cap(E+\e C)|-|E|)\,.
\end{equation}
Actually, this definition is not very practical, since it can change
drastically with Lebesgue-negligible changes of the set $E$. For
this reason, we introduce the functional, defined for
a measurable function $u$:
\begin{equation}\label{defFe}
F_{\e, C}(u;\Om)\ :=\ \frac{1}{\e}\int_{\Om}\big( \esssup_{\Om\cap (x-\e C)}u-u(x)\big)\,dx
\end{equation}
which takes values in $[0,+\infty]$. Notice that one can check easily
(using Fatou's lemma)
that $F_{\e,C}$ is l.s.c.~in $L^1_\loc(\Om)$. We then define
\begin{equation}\label{defsmec}
\SMeC(E;\Om)\ :=\ F_{\e, C}(\chi_E;\Om)\,.
\end{equation}
It is also easy to check that the definition coincides with
$\SMeC^0$ on smooth sets, and in general for a measurable set $E$
we have
\[
\SMeC(E;\Om)\ =\ \min_{|E'\triangle E|=0} \SMeC^0(E';\Om)\ =\ \SMeC^0(E^1;\Om)
\ =\ \SMeC^0(\Om\setminus E^0;\Om)
\]
where $E^1$ (resp., $E^0$) is the set of points of Lebesgue
density $1$ (resp., $0$) in $\Om$ (see~\eqref{defEt}),
and $A\triangle B$ denotes the symmetric difference
$(A\setminus B)\cup (B\setminus A)$ of the sets $A$
and $B$. Eventually, one can check easily that
$F_{\e, C}$ satisfies a \textit{generalized coarea formula} \cite{Visintin,CGL}:
for any function $u\in L^1_\loc(\Om)$,
\begin{equation}\label{gcoarea}
F_{\e,C}(u;\Om)\ =\ \int_{-\infty}^\infty \SMeC(\{u>s\};\Om)\,ds\,.
\end{equation}
Indeed, for any integer $k$ and a.e.~$x\in\Om$ with $u(x)>-k$,
\[
\left(\esssup_{\Om\cap(x-\e C)} u\right) + k
 \ =\ \int_{-k}^{+\infty} \esssup_{\Om\cap(x-\e C)}
\chi_{\{u > s\}}\,ds
\]
so that
\[
\esssup_{\Om\cap(x-\e C)} u-u(x)\ =\ \int_{-k}^{+\infty}
\esssup_{\Om\cap(x-\e C)} \chi_{\{u > s\}} - \chi_{\{u>s\}}(x)\,ds\,,
\]
and sending $k\to+\infty$ we deduce~\eqref{gcoarea}.

Eventually, one can show that for any measurable sets, $E$ and $F$,
\[
\SMeC(E\cup F;\Om)+\SMeC(E\cap F;\Om)\ \le\ \SMeC(E;\Om)+\SMeC(F;\Om)
\]
from which it follows that $F_\e$ is convex on $L^1_\loc(\Om)$, see~\cite{CGL}
for details.

Before stating our results we recall the following definition:
we say that a family of functionals $(G_\e)_{\e>0}$
defined on the Lebesgue measurable subsets of $\R^n$ $\Gamma$-converges
to $G$ in $L^1_\loc(\Om)$ as $\e\to 0$ if
\begin{itemize}
\item  for any Lebesgue measurable set $E$, and for any family $(E_\e)_{\e>0}$
of Lebesgue measurable sets $E_\e\to E$ in $L^1_\loc(\Om)$ as $\e\to 0$,
we have
\begin{equation*}\label{}
G(E)\ \le\ \liminf_{\e\to 0} G_\e(E_\e)\,,
\end{equation*}
\item for any Lebesgue measurable set $E$, there exists a family of Lebesgue measurable sets $(E_\e)_{\e>0}$
such that $E_\e\to E$ in $L^1_\loc(\Om)$  as $\e\to 0$ and
\begin{equation*}\label{}
\limsup_{\e\to 0} G_\e(E_\e)\ \le\ G(E).
\end{equation*}
\end{itemize}

In the same way, we say that a family of functionals $(F_\e)_{\e>0}$ defined on
$L^1_\loc(\Om)$ $\Gamma$-converges to $F$ in $L^1_\loc(\Om)$  as $\e\to 0$ if
\begin{itemize}
\item  for any $u\in L^1_\loc(\Om)$, and for any family $(u_\e)_{\e>0}$
of elements of $L^1_\loc(\Om)$ such that $u_\e\to u$ in $L^1_\loc(\Om)$
as $\e\to 0$, we have
\begin{equation*}\label{}
F(u)\ \le\ \liminf_{\e\to 0} F_\e(u_\e)\,,
\end{equation*}
\item for any  $u\in L^1_\loc(\Om)$, there exists a family $(u_\e)_{\e>0}$
of elements of $L^1_\loc(\Om)$ such that $u_\e\to u$ in $L^1_\loc(\Om)$
as $\e\to 0$ and
\begin{equation*}\label{}
\limsup_{\e\to 0} F_\e(u_\e)\ \le\ F(u).
\end{equation*}
\end{itemize}

We will show the following result:
\begin{theorem}\label{ThGCV}
As $\e\to 0$, $\SMeC$ and $\SMeC^0$ $\Gamma$--converge to
\begin{equation}\label{defPhC}
\PhC(E;\Omega) \ :=\ \begin{cases}
\displaystyle \int_{\partial^* E\cap\Om} h_C(\nu_E(x))\,d\hs^{n-1}(x)
& \textrm{ if $E$ has finite perimeter in $\Om$}\,,
\\ +\infty &\textrm{ else}
\end{cases}
\end{equation}
in $L^1_\loc(\Om)$, where $\nu_E(x)$ is the
\textit{outer}\footnote{Observe that with this classical
but not so natural choice, we have $\PhC(E)=\int_\Om h_C(-D\chi_E)$.}
normal to $\partial^* E$ at $x$.
Moreover, if $\{E_\e\}_{\e>0}$ are sets with locally finite measure and
$\sup_{\e>0}\SMeC(E_\e;\Om)<+\infty$, then, up to subsequences, $E^1_\e$
converge to some set $E$ in $L^1_\loc(\Om)$.
\end{theorem}
In particular, we deduce from~\cite[Prop.~3.5]{CGL}:
\begin{corollary}\label{CorGCV}
As $\e\to 0$, $F_{\e,C}$ $\Gamma$--converges to
\begin{equation*}
TV_{-C}(u;\Om) \ :=\ \begin{cases}
\displaystyle \int_{\Om} h_C(-Du)& \textrm{ if }u\in BV(\Om)\,,
\\ +\infty &\textrm{ else}
\end{cases}
\end{equation*}
in $L^1_\loc(\Om)$, where $h_C(-Du)$ stands for $h_C\big(-\frac{dDu}{d|Du|}\big)\,d|Du|$~\cite{DemengelTemam}. Moreover, if $\{u_\e\}_{\e>0}$ are functions in
$L^1_\loc(\Om)$ with $\sup_{\e}F_{\e,C}(u_\e;\Om)<+\infty$, then $\{u_\e\}_{\e>0}$
is precompact in $L^1_\loc(\Om)$.
\end{corollary}
Indeed, Proposition~3.5 in~\cite{CGL} shows that for functionals
such as $F_{\e,C}$ which satisfy~\eqref{gcoarea},
the $\Gamma$-convergence of the
functionals restricted to characteristics functions
of sets to some limit (Theorem~\ref{ThGCV}) is sufficient to
imply the $\Gamma$-convergence of the full
functionals to a limit which is precisely the extension by co-area formula
of the previous limit (here $TV_{-C}$ is the extension of~\eqref{defPhC}).
The compactness is shown in Section~\ref{secproofgcv}.

For any measurable set $E$ we can also consider
$$
\mathcal M_{\e,C}(E;\Om)\ := \ (\SMeC(E;\Om)+\SMeC(\Om\setminus E;\Om))/2.
$$
From Theorem \ref{ThGCV} the following Corollary follows easily: 

\begin{corollary}\label{cor}
As $\e\to 0$, $\mathcal M_{\e,C}$ $\Gamma$--converges to $(\PhC(E)+\PhC(\Om\setminus E))/2$ in $L^1_\loc(\Om)$.
\end{corollary}

Concerning the pointwise convergence of $\SMeC^0$, we also have
the following interesting result, from which the $\Gamma$-$\limsup$
inequality in Theorem~\ref{ThGCV} follows in a straightforward way.
\begin{theorem}\label{ThPW}
Assume that the set $E$ is a finite-perimeter set such that
\begin{equation}\label{pwlim}
\lim_{\e \to 0} \SM_{\e,B(0,1)}^0(E;\Om)\ =\ \P(E;\Om)\,.
\end{equation}
Then,
\begin{equation}\label{pwlimh}
\lim_{\e \to 0} \SMeC^0(E;\Om)\ =\ \PhC(E;\Om)\,.
\end{equation}
\end{theorem}
\begin{oss}
If we assume moreover that $C$ is $C^{1,1}$ and ``elliptic'' (precisely:
that $h_C^2/2$ is both uniformly convex and with Lipschitz gradient),
then the two assertions in Theorem \ref{ThPW} should in fact be
equivalent. To check this, it requires to adapt the proof by replacing
the Euclidean ball and distance with $C$ and the corresponding
anisotropic (nonsymmetric) distance. The smoothness and ellipticity
are required because, in the proof, we use the ellipticity of the
distance to control the difference of measure between a flat surface,
orthogonal to a given vector $\nu$, and a slanted $C^{1,1}$ surface
with at least a point with normal $\bar\nu\neq \nu$, see Fig.~\ref{FigBall}
and the proof of~\eqref{estimbelowFs} for details.
\end{oss}
\begin{oss}\label{remMinkIso}
The sets which satisfy~\eqref{pwlim} are studied in~\cite{ACV}.
A sufficient condition is that the Minkowski content of
the reduced boundary coincides with its $(n-1)$--dimensional
measure, that is,
\begin{equation}\label{minkok}
\lim_{\e\to 0}\frac{|\{x\in \Om\,:\, \dist(x,\partial^* E)\le \e\}|}{2\e}
\ =\ \P (E;\Om)\,.
\end{equation}
The proof is quite elementary (see~\cite[Thm~5]{ACV}).
Indeed, we observe first that~\eqref{minkok} implies that
$\{\dist(\cdot,\partial^* E)=0\}$
has zero Lebesgue measure. Then, we can introduce for $x\in \Om$
the (essential) signed distance function
\[
d_E(x)\,:=\,\dist(x,E^1)-\dist(x,E^0)
\]
to the boundary. If $\dist(x,E^1)>0$ then clearly $x$ is in the
interior of $E^0$, in the same way if $\dist(x,E^0)>0$,
$x$ is in the interior of $E^1$.
Hence $|d_E(x)|=\dist(x,E^1)+\dist(x,E^0)$, and in particular
$d_E(x)=0$ if and only if for any $\rho>0$,
$B(x,\rho)$ contains points of
both $E^0$ and $E^1$,
which in turn is equivalent to $\dist(x,\partial^* E)=0$. It follows
that $|d_E(x)|=\dist(x,\partial^* E)$ (and, by~\eqref{minkok} $\{d_E=0\}$ is
negligible, showing also that $d_E$ is the classical signed distance function
to $\partial E^1$ in $\Om$).
Thanks to the co-area formula and the fact that $|\nabla d_E|=1$ a.e.,
we have both
\begin{equation*}
\begin{array}{ll}
\displaystyle \frac{1}{\e}|\{x\in \Om\,:\, 0< d_E(x)< \e\}|
\ =\ \int_0^1 \P (\{d_E<\e s\}; \Om)\,ds\,,\\
\displaystyle \frac{1}{\e}|\{x\in \Om\,:\, -\e< d_E(x)< 0\}|
\ =\ \int_{-1}^0 \P (\{d_E<\e s\}; \Om)\,ds\,.
\end{array}
\end{equation*}
Now for any $s$, $\{d_E<\e s\}$ goes to either $\{d_E<0\}$ (if $s<0$)
or $\{d_E\le 0\}$ (if $s>0$) as $\e\to 0$, which both
coincide to $E$ up to a negligible
set. Hence by Fatou's lemma,
\begin{equation}\label{limminkpm}
\begin{aligned}
&\P(E;\Om)
\ \le\ \liminf_{\e\to 0} \frac{1}{\e}|\{x\in \Om\,:\, 0 <d_E(x)<\e\}|,\\
&\P(E;\Om)
\ \le\ \liminf_{\e\to 0} \frac{1}{\e}|\{x\in \Om\,:\, 0 <-d_E(x)< \e\}|.
\end{aligned}
\end{equation}
Together with~\eqref{minkok}, it follows that the inequalities
in~\eqref{limminkpm} must in fact be equalities, and
the $\liminf$s are limits. In particular, we deduce~\eqref{pwlim}.
\end{oss}
If $\partial E$ is compact and rectifiable (that is, included in the image of a Lipschitz map from $\R^{n-1}$ to $\R^n$), and $\hs^{n-1}(\partial E\setminus \partial^* E)=0$, then the Minkowski content coincides with the perimeter, see \cite[Thm~3.2.39 p.~275]{Federer}, and the previous analysis applies.
It is easy to build examples, though, where this is not true and still,
\eqref{pwlimh} holds, see again~\cite{ACV}.
\begin{oss}\label{remConv}
In case $E$ is a convex body, then it is well known that (see~\cite{Schneider})
\[
|E+\e C|=|E|+\e \PhC(E;\Om)+O(\e^2)\,.
\]
\end{oss}

As before, for any measurable set $E \subset \Om$ we let
\[
\mathcal M_{\e,C}^0(E;\Om):=(\mathcal {SM}_{\e,C}^0(E;\Om)+\mathcal {SM}_{\e,C}^0(\Om\setminus E;\Om))/2.
\]
Then the following pointwise convergence result holds.

\begin{theorem}\label{ThPW1}
Assume that the set $E$ is a finite-perimeter set such that
\[
\lim_{\e \to 0} \mathcal M_{\e,B(0,1)}^0(E;\Om)\ =\ \P(E;\Om)\,.
\]
Then,
\[
\lim_{\e \to 0} \mathcal M_{\e,C}^0(E;\Om)\ =\ (\PhC(E;\Om)+\PhC(\Om\setminus E;\Om))/2\,.
\]
In particular, we get
\[
\lim_{\e \to 0} \mathcal M_{\e,C}^0(E;\Om)=\int_{\partial^* E} \frac{h_C(\nu_E(x))+h_C(-\nu_E(x))}{2}\,d\hs^{n-1}(x).
\]
\end{theorem}

\section{Proof of the results}\label{secProofs}

\subsection{Proof of Theorem~\ref{ThGCV}}\label{secproofgcv}


In order to prove the $\Gamma$--convergence, we must show that for
any $E$,
\begin{itemize}
\item  if $E_\e\to E$ in $L^1_\loc(\Om)$,
then
\begin{equation}\label{Glinf}
\PhC(E;\Om)\ \le\ \liminf_{\e\to 0} \SMeC(E_\e;\Om)\,,
\end{equation}
\item and that there exists $E_\e\to E$ with
\begin{equation}\label{Glsup}
\limsup_{\e\to 0} \SMeC(E_\e;\Om)\ \le\ \PhC(E;\Om).
\end{equation}
\end{itemize}

As it is standard that one can approximate any set $E$ with finite perimeter by means of smooth (enough)
sets such that $\P(E_k;\Om)\to \P(E;\Om)$ (for instance, minimizers
of $\min_F\P(F;\Om)+k|E\triangle F|$ will have a $C^1$ boundary, up to
a compact singular set of small dimension)
then \eqref{Glsup} will follow, using a diagonal argument
and Remark~\ref{remMinkIso}, from Theorem~\ref{ThPW} (which we will prove later on).

Hence, we focus on the proof of~\eqref{Glinf}. We will also prove,
simultaneously, the last claim of the theorem, which is the compactness
of a family $(E_\e)$ with equibounded energies.
Let us introduce
the anisotropic (essential) distance function to a set $E$:
\[
\distC(x,E)\ :=\ \essinf_{y\in E}\,h^\circ_C(x-y)\,.
\]
(Equivalently, this is the $h^\circ_C$-distance to the set $E^1$ of
points where the Lebesgue density of $E$ is $1$, or to the complement
of $E^0$.)
Then, $\distC(x,E)< \e$ if and only if there
exists a set of positive measure in $E$ of points $y$ with
$h^\circ_C(x-y)<\e$, or, in other words, such that $x-y\in \e C$,
which is equivalent to say that $x\in E^1+\e C$. In particular,
it follows that
\[
(E^1+\e C\setminus E^1)\cap \Om=\{x\in\Om\setminus E^1\,:\,
\distC(x,E)< \e\}.
\]
On the other hand, if one lets $d(x):=\distC(x,E)$, it is
standard that $d$ is Lipschitz and that $h_C(\nabla d)=1$
a.e.~in $\{ d>0\}$, and $h_C(\nabla d)=0$ a.e.~in $\{d=0\}\supset E^1$.
The proof of this fact follows the same lines
as in~\cite{AmbrosioDistance}: first, for any
$x,y\in\Om$, if $\delta>0$, one can find a set with positive measure
in $E$ of points $y'$ with $d(y)\le h^\circ_C(y-y') \le d(y)+\delta$.
Then, for these points,
\begin{equation}\label{dlip}
d(x)-d(y) \ \le\ h^\circ_C(x-y')-h^\circ_C(y-y')+\delta
\ \le\ h^\circ_C(x-y)+\delta
\end{equation}
and sending $\delta$ to zero and using~\eqref{comparh}, it follows that $d$
is Lipschitz. Moreover, $\nabla d=0$ a.e. in~$\{d=0\}$.
Now, from~\eqref{dlip} we also see that $d(x+tz)-d(x)\le th^\circ_C(z)$
for all $z$; therefore, if $d$ is differentiable at $x$ it follows that
$\nabla d(x)\cdot z\le 1$ for all $z\in C$, hence $h_C(\nabla d(x))\le 1$.

We show the reverse inequality for points $x$ where $d(x)>0$:
for such a point, there exists $y\in \ov{E^1}$ with $d(x)=h^\circ_C(x-y)$.
For each $x'\in (y,x]$ (which means that $x' \neq y$ and $x'$ lies on the line segment with extreme points $y$ and $x$), one has
$d(x')=h^\circ_C(x'-y)>0$, otherwise there would exist $y'$ with
$h^\circ_C(x'-y')<h^\circ_C(x'-y)$, but then, it would follow
that
$$
h^\circ_C(x-y')\le h^\circ_C(x-x')+h^\circ_C(x'-y')<h^\circ_C(x-x')+h^\circ_C(x'-y)
=h^\circ_C(x-y)
$$
since $x'\in (y,x]$, a contradiction.
It follows that for $z=x-y$, $t\in (0,1)$,
$$
d(x-tz)=h^\circ_C(x-tz-y)=(1-t)d(x),
$$
and if
in addition $x$ is a point of differentiability,
it follows that
$$
-\nabla d(x)\cdot z=-d(x)=-h^\circ_C(z).
$$
But
since $h_C(\nabla d(x))\le 1$ and $z/h^\circ_C(z)\in C$, it
follows that $h_C(\nabla d(x))= 1$. If moreover $h^\circ_C$ is differentiable
as well in  $x-y$, we find in addition that
$\nabla d(x)=\nabla h^\circ_C(x-y)$. If $h_C$ is differentiable
in $\nabla d(x)$, we find that $y=x-d(x)\nabla h_C(\nabla d(x))$
and in particular, in that case, the projection $y$ must be unique.
For a general convex set $C$ this might not be the case, even at
points of differentiability.

Let us now show~\eqref{Glinf} and the compactness.
We let $\{E_\e\}_{\e>0}$ be a family of sets, with
$\liminf_{\e\to0} \SMeC(E_\e;\Om)<+\infty$.
We consider a subsequence $E_k:=E_{\e_k}$ such that this
$\liminf$ is in fact a limit. We will show both that, up to subsequences,
it converges to a set $E$ in $L^1_\loc(\Om)$ and that~\eqref{Glinf} holds.
We have
\begin{equation*}
((E_k^1+\e_k C)\setminus E_k^1)\cap \Om\ \supseteq\ \{x\in\Om\,:\,
0<\distC(x,E_k)<\e_k\}
\end{equation*}
(the difference being the possible set of points $x\not\in E_k^1$
with $\distC(x,E_k)=0$).
It follows, letting $d_k(x):= \min\{\distC(x,E_k)/\e_k,1\}$,
\begin{equation*}
\SM_{\e_k,C}(E_k;\Om)\ \ge\ \frac{1}{\e_k}
\int_{\{0<\distC(\cdot,E_k)<\e_k\}} h_C(\nabla \distC(x,E_k))\,dx
\ =\ \int_\Om  h_C(\nabla d_k(x))\,dx\,.
\end{equation*}
In particular, $(d_k)_{k\ge 1}$ have equibounded total variation:
we may assume that a subsequence (not relabelled)
converges to some limit $d$, with values in $[0,1]$, in $L^1_\loc(\Om)$.
(And, in fact, we may even assume that the convergence is pointwise,
out of a negligible set.)

By assumption, $|\{0<d_k<1\}|\le |\{d_k<1\}\setminus E^1_k|\le c\e_k$,
in particular we deduce easily
that $d\in \{0,1\}$ a.e.~in $\Om$ (for instance, by checking
that $d_k(1-d_k)\to 0$). We call $E=\{d=0\}$. In particular,
$\chi_E=1-d$. Observe that
if $B$ is a ball in $\R^n$,
\begin{multline*}
\int_B |\chi_{E_k}-\chi_E|\,dx\\
=\,\int_{B\cap E^1_k} |d_k-d|\,dx
\,+\, \int_{B\cap (\{d_k<1\}\setminus E^1_k)} |\chi_E|
\,+\, \int_{B\cap \{d_k=1\}} |d_k-d|\,dx
\\ \le\ \|d_k-d\|_{L^1(B\cap\Om)} \,+\,c\e_k\ \to\ 0
\end{multline*}
as $k\to\infty$, so that $E_k\to E$ in $L^1_\loc(B)$, hence showing the
compactness.

Thanks to Reshetniak's lower semicontinuity Theorem, it follows
from the $L^1_\loc$-convergence of $d_k$ to $1-\chi_E$ that
\begin{equation*}
\int_\Om h_C(-D\chi_E)\ \le\ \liminf_{k\to\infty} \int_\Om h_C(\nabla d_k(x))\,dx\,
\ \le\ \lim_{k\to\infty} \SM_{\e_k,C}(E_k;\Om)\,.
\end{equation*}
Since $\int_\Om h_C(-D\chi_E)=\PhC(E;\Omega)$, \eqref{Glinf} follows.
\smallskip

To extend the compactness result to Corollary~\ref{CorGCV}, one can
consider for each $\delta>0$ and function $u_\e$ a function
$$
u_\e^\delta=\sum_{k\in\Z} s_k\chi_{\{s_{k+1}\ge u_\e>s_k\}},
$$
where $s_k\in (k\delta,
(k+1)\delta)$ is a level  appropriately chosen so that
$$
\SMeC(\{u_\e>s_k\};\Om)\le (1+\sup_{\e>0} F_{\e,C}(u_\e;\Om))/\delta.
$$
Then,
the previous compactness result (and a diagonal argument)
shows that there exists a a positive infinitesimal sequence $\e_k$ such that $u_{\e_k}^{1/n}$
converges to some $u^{1/n}$ in $L^1_\loc(\Om)$, for all $n\ge 1$.
Since $\|u_{\e_k}^{1/m}-u_{\e_k}^{1/n}\|_\infty\le 2/\min\{m,n\}$
and $\|u^{1/m}-u^{1/n}\|_\infty\le 2/\min\{m,n\}$ for all $m,n,k$,
we easily deduce that (up to a subsequence), there exists $u$
such that $u_{\e_k}\to u$ in $L^1_\loc(\Om)$.
\bigskip

As already mentioned, the proof of \eqref{Glsup} will follow from
Theorem~\ref{ThPW}, which is proved in the next Section.

\subsection{Proof of Theorem~\ref{ThPW}}

Now, we consider a set $E\subset \Om$ such that \eqref{pwlim} holds.
We will identify $E$ with the set of points where its Lebesgue density
is $1$, moreover, a necessary condition for~\eqref{pwlim} is that
$\ov{E}=\bigcap_{\e>0} E+B(0,\e)$ coincides with $E$ up to a negligible
set, in other words, $|\ov{E}\setminus E|=0$.

A first remark is that, clearly, using \eqref{comparh},
$$a\SM^0_{a\e, B(0,1)}(E;\Om)\le \SMeC^0(E;\Om)\le b\SM^0_{b\e, B(0,1)}(E;\Om)$$
hence any limit of $\SMeC^0(E;\Om)$ is in between $a\P(E;\Om)$ and
$b\P(E;\Om)$. In particular, we can introduce the non-negative
measures
\begin{equation*}
\mu_\e\ :=\ \frac{1}{\e}\left( \chi_{E+\e C}-\chi_E\right)\hs^n
\end{equation*}
which are equibounded, since by definition $\mu_\e(\Om)=\SMeC^0(E;\Om)$.
Then, up to a subsequence, we have $\mu_{\e_k}\wtos \mu$
as measures in $\Om$, with $a\hs^{n-1}\mrestr \partial^* E\le
\mu\le b\hs^{n-1}\mrestr \partial^* E$.
In order to prove the result, we need to show that $\mu$
is equal to $h_C(\nu_E)\hs^{n-1}\mrestr\partial^* E$.

For this purpose, we introduce the Besicovitch
derivative $g$
of the measure $\mu$ w.r.t.~$\hs^{n-1}\mrestr\partial^* E$, defined by
\begin{equation*}
g(x)\ =\ \lim_{\rho\to 0} \frac{\mu(B(x,\rho))}{\hs^{n-1}(\partial^*E\cap
B(x,\rho))}
\end{equation*}
(observe that  $g(x)\in [a,b]$)
and which is defined for $\hs^{n-1}$--a.e.~$x\in \Om\cap\partial^* E$. Moreover,
since $\partial^* E$ is rectifiable, it is also given a.e.~by
\begin{equation}\label{RadonNyk}
g(x)\ =\ \lim_{\rho\to 0} \frac{\mu(B(x,\rho))}{\alpha_n \rho^{n-1}}
\end{equation}
where $\alpha_n=\omega_{n-1}$ is the measure of the $n-1$ dimensional unit ball.

As mentioned,
Theorem~\ref{ThPW} will follow if we can show that $g(x)= h_C(\nu_E(x))$
for $\hs^{n-1}$--a.e.~$x\in \partial^* E$.
Observe that from~\eqref{Glinf}, it
follows that $g(x)\ge h_C(\nu_E(x))$ for $\hs^{n-1}$-a.e.~$x\in\partial^* E$,
so that we just need to show that $g(x)\le h_C(\nu_E(x))$ for $\hs^{n-1}$-a.e.~$x\in\partial^* E$.

A first step, which is classical, is to blow-up the space around
a point $\bar x$ where the Besicovitch
derivative exists, and by a diagonal argument, to consider the
situation where the set is close to a half space,
orthogonal to $\nu_E(\bar x)$.
We thus fix from now on a point
$\bar x\in \partial^* E$ where~\eqref{RadonNyk} holds.

We recall that
\begin{equation*}
 \lim_{\rho\to 0} \frac{\hs^{n-1}(\partial^*E\cap
B(\bar x,\rho))}{\alpha_n \rho^{n-1}} \ =\ 1
\end{equation*}
and
\begin{equation}\label{bupbup}
 \lim_{\rho\to 0} \int_{B(0,1)} \left|\chi_{\frac{E-\bar x}{\rho}}
 - \chi_{\{y\,:\,y\cdot\nu_E(\bar x)\le 0\}} \right|\,dy\ =\ 0
\end{equation}
hold. We denote $\nu=\nu_E(\bar x)$ and without loss of generality
we will assume that it is the direction of the last coordinate $x_n$.
We will use the notation $x=(x',x_n)\in \R^{n-1}\times \R$ to distinguish
between the component $x'\perp \nu$ and $x_n$ (along $\nu$) of a point
$x\in\R^n$.

We also introduce the measures
\begin{equation*}
\lambda_\e\ :=\ \frac{1}{\e}\left( \chi_{E+\e B(0,b)}-\chi_E\right)\hs^n \ \ge\ \mu_\e\,,
\end{equation*}
the main assumption of Theorem~\ref{ThPW} ensures that these
measures converge weakly-$*$ to $\lambda=b\hs^{n-1}\mrestr \partial^* E$
as $\e\to 0$. Now we use a classical procedure: since for a.e.~$\rho>0$,
\[
\mu_{\e_k}(B(\bar x,\rho))\to\mu(B(\bar x,\rho)) \quad \textrm{and} \quad \lambda_{\e_k}(B(\bar x,\rho))\to\lambda(B(\bar x,\rho)),
\]
(see \cite{AFP-libro} Proposition 1.62 and Example 1.63)
we can build an infinitesimal sequence
$(\rho_k)_{k\in\N}$ with $\e'_k=\e_k/\rho_k\to 0$ as $k\to \infty$ such that
\begin{equation}\label{derivC}
\lim_{k\to\infty} \frac{\mu_{\e_k}(B(\bar x,\rho_k))}{\alpha_n \rho_k^{n-1}}
 \ =\ g(\bar x)
\end{equation}
and
\begin{equation}\label{derivbB}
\lim_{k\to\infty} \frac{\lambda_{\e_k}(B(\bar x,\rho_k))}{\alpha_n \rho_k^{n-1}}
 \ =\ b.
\end{equation}
The rest of the proof would be relatively easy
if we could ensure that $\e_k\sim \rho_k$ as $k\to\infty$, using
then a blow-up argument.
The reason is that in this case, at the scale $\rho_k$, the set
$E$ would look like a half-space while $(E+bB(0,\e_k))\setminus E$ would look
like a strip, of constant width $\sim b(\e_k/\rho_k)$.
The fact that the volume of this strip goes precisely to the
volume of a straight strip (which \eqref{derivbB} tells us)
would then imply that it is essentially
straight, up to a small error. This would, in particular, show that
at the scale $\e_k$, $\partial E$ is almost flat
and we would be able
to estimate precisely the volume of $(E+\e_kC)\setminus E$.

However, this is not clear in general, and
we need to consider the general situation, where $\e_k=o(\rho_k)$, hence
$\e'_k\to 0$.
The workaround will be to consider, after a blow-up at scale $\rho_k$,
a covering of the (flat) limit surface with cubic regions of scale $\e'_k$
and show that ``most of'' these regions are good, meaning that they
can be roughly analyzed at scale $\e'_k$
with the arguments previously mentioned,
while the other regions are not enough
to contribute significantly to the limit.

As is usual, we do a blow-up using the change of
variables $x=\bar x+\rho_k y$.
We let $E_k=(E-\bar x)/\rho_k$, and we observe that from \eqref{derivC},
\eqref{derivbB}~and \eqref{bupbup},
\begin{eqnarray}
\label{eq:BUC}
\lim_{k\to\infty} \frac{1}{\alpha_n \e'_k}\int_{B(0,1)} (\chi_{E_k+\e'_k C}-\chi_{E_k})\,dy
 \ =\ g(\bar x), & \\[1mm]
\label{eq:BUb}
\lim_{k\to\infty} \frac{1}{\alpha_n \e'_k}\int_{B(0,1)} (\chi_{E_k+\e'_k B(0,b)}-\chi_{E_k})\,dy
 \ =\ b, & \\[1mm]
\label{eq:BUL1}
 \lim_{k\to\infty}
\int_{B(0,1)} \left|\chi_{E_k}-\chi_{\{y\,:\,y\cdot\nu\le 0\}} \right|\,dy\ =\ 0 .&
\end{eqnarray}
Moreover, for any $\beta>0$ (small), one can check easily that if
we replace in~\eqref{eq:BUb} and~\eqref{eq:BUL1} $B(0,1)$ with
$B(0,1-\beta)$, or even with $C(0,1-\beta):=\{(y',y_n)\in B(0,1)\,:\,
|y'|\le 1-\beta\}$, \eqref{eq:BUL1} still holds and
the right-hand side in~\eqref{eq:BUb} is replaced with $b(1-\beta)^{n-1}$.
Indeed, it follows from~\eqref{Glinf} (with $C=B(0,b)$) and~\eqref{eq:BUL1} that
for any open set $A\subseteq B(0,1)$,
\begin{equation*}
b\hs^{n-1}(A\cap \{y\cdot\nu=0\})\ \le\ \liminf_{k\to\infty}
 \frac{1}{\e'_k}\int_{A} (\chi_{E_k+\e'_k B(0,b)}-\chi_{E_k})\,dy\,.
\end{equation*}
Together with~\eqref{eq:BUb}, we deduce that
as soon as $\hs^{n-1}(\partial A\cap \{y\cdot\nu = 0\})=0$,
\begin{equation}\label{eq:BUblocal}
b\hs^{n-1}(A\cap \{y\cdot\nu=0\})\ =\ \lim_{k\to\infty}
 \frac{1}{\e'_k}\int_{A} (\chi_{E_k+\e'_k B(0,b)}-\chi_{E_k})\,dy\,.
\end{equation}

We fix a (small) value of $\beta>0$.
Then, we choose a value
$\theta>10b$ and we consider the points $z\in \Z^{n-1}$ such
that the hypersquares $(\theta\e'_k (z+(0,1)^{n-1}))\times\{0\}$
are contained in $B(0,1-\beta)$. There is a finite number $N_k$
of such squares and we enumerate the corresponding
points $\{z^k_1,\dots,z^k_{N_k}\}$. For $i=1,\dots,N_k$, we let
\begin{equation*}
C_i^k\ =\ [(\theta\e'_k( z^k_i+(0,1)^{n-1}))\times\R]\cap B(0,1)\,,
\qquad {C'_i}^k\ =\ (\theta\e'_k( z^k_i+(0,1)^{n-1}))\times\{0\}\,.
\end{equation*}
We then let
\begin{equation}\label{defaik}
a_i^k\ =\ \int_{C_i^k} \left|\chi_{E_k}-\chi_{\{y_n\le 0\}} \right|\,dy
\ \le\ 2(\theta \e'_k)^{n-1}
\end{equation}
and $\delta_k=\sum_{i=1}^{N_k} a_i^k$: from~\eqref{eq:BUL1} we know that
$\delta_k\to 0$ as $k\to\infty$. We then consider
$$
Z_k=\{i =1,\dots,N_k\,:\, a_i^k \le \sqrt{\delta_k}(\e'_k)^{n-1}\},
$$
$Z'_k=\{1,\dots,N_k\}\setminus Z_k$.
It follows that $\delta_k\ \ge\ \sqrt{\delta_k}(\e'_k)^{n-1}\# Z'_k$ and then
\begin{equation}\label{controlZpk}
(\e'_k)^{n-1}\# Z'_k\ \le\ \sqrt{\delta_k}
\end{equation}
which gives a control on the ``bad'' surface, of the cylinders
$C^k_i$ where the integral $a_i^k$ is ``large''.
On all the other cylinders, if we blow-up the coordinates at
scale $\e'_k$ we will still have that $E_k$ is close,
in some sense, to $\{y_n\le 0\}$.

For each $i=1,\dots,N_k$, 
we have
\begin{multline}\label{controlset}
\frac{1}{\e'_k}\int_{C_i^k}\left|\chi_{E_k+\e'_k B(0,b)}-\chi_{E_k}\right|\,dy
\,=\, \frac{1}{\e'_k}\int_{{C'_i}^k} \int_{-\sqrt{1-{y'}^2}}^{\sqrt{1-{y'}^2}}
\left|\chi_{E_k+\e'_kB(0,b)}-\chi_{E_k}\right|\,dy_n dy'
\\ \ge\, b\left|
\left\{y'\in {C'_i}^k\,:\, \hs^1((\{y'\}\times\R)\cap(B(0,1)\cap E_k))>0,
\hs^1((\{y'\}\times\R)\cap(B(0,1)\setminus E_k))\ge \e'_k b \right\}
\right|
\end{multline}
since clearly, each time a point $(y',y_n)\in E_k$, then $(y',y_n+s)\in
E_k+\e'_kB(0,b)$ for $| s|\le \e'_kb$. We denote by $D_i^k$ the set
in the right-hand side of~\eqref{controlset}. For $y'\not\in D_i^k$,
\begin{equation*}
\int_{-\sqrt{1-{y'}^2}}^{\sqrt{1-{y'}^2}}
\left|\chi_{E_k}-\chi_{\{y_n\le 0\}} \right|\,dy\ \ge\ \sqrt{\beta}-\e'_k b
\ \ge\ \frac{\sqrt{\beta}}{2}
\end{equation*}
as soon as $\e'_k\le \sqrt{\beta}/(2b)$ (which we assume in the sequel).
It follows that $|{C'_i}^k\setminus D_i^k|\le 2 a_i^k/\sqrt{\beta}$, hence
if $i\in Z_k$, so that $a_i^k\le  \sqrt{\delta_k}(\e'_k)^{n-1}
=\sqrt{\delta_k}|{C'_i}^k|/\theta^{n-1}$, we get that
\begin{equation}\label{controlsetlwbd}
b|D_i^k|\ \ge\ b|{C'_i}^k|\left(1- \frac{2\sqrt{\delta_k}}{\theta^{n-1}\sqrt{\beta}}
\right)\ =\ b|{C'_i}^k|\left(1- K\sqrt{\delta_k}\right)\,.
\end{equation}
To sum up, \eqref{controlset} and~\eqref{controlsetlwbd} show that for
any $i\in Z_k$,
\begin{equation*}
\frac{1}{\e'_k}\int_{C_i^k}\left|\chi_{E_k+\e'_k B(0,b)}-\chi_{E_k}\right|\,dy
\ \ge\ b|{C'_i}^k|\left(1- K\sqrt{\delta_k}\right)\,.
\end{equation*}
In particular, it follows that (using~\eqref{controlZpk})
\begin{multline}\label{bigestimk}
\liminf_{k\to\infty} \sum_{i\in Z_k}
\frac{1}{\e'_k}\int_{C_i^k}\left|\chi_{E_k+\e'_k B(0,b)}-\chi_{E_k}\right|\,dy
\\ \ge\ b\lim_{k\to\infty}
\left(\left|\bigcup_{i=1}^{N_k}{C'_i}^k\right|-(\theta\e'_k)^{n-1}\# Z'_k\right)
\left(1- K\sqrt{\delta_k}\right)\ =\ \alpha_n b(1-\beta)^{n-1}
\end{multline}
and together with~\eqref{eq:BUblocal}
(which bounds the $\limsup$, for $A=C(0,1-\beta)$) we deduce that
\begin{equation*}
\delta'_k\ :=\
 \sum_{i\in Z_k}
\left(
\frac{1}{\e'_k}\int_{C_i^k}\left|\chi_{E_k+\e'_k B(0,b)}-\chi_{E_k}\right|\,dy
- b|{C'_i}^k|\left(1- K\sqrt{\delta_k}\right)
\right)\ \stackrel{k\to\infty}{\longrightarrow} \ 0.
\end{equation*}

So now we introduce
\begin{equation*}
\tilde{Z}_k\,=\,\left\{i \in Z_k\,:\,
\frac{1}{\e'_k}\int_{C_i^k}\left|\chi_{E_k+\e'_k B(0,b)}-\chi_{E_k}\right|\,dy
- b|{C'_i}^k|\left(1- K\sqrt{\delta_k}\right)\,\le\,
(\e'_k)^{n-1}\sqrt{\delta'_k}
\right\},
\end{equation*}
and its complement $\tilde{Z}'_k=Z_k\setminus \tilde{Z_k}$. Then
as before, one sees that $\delta'_k\ \ge\ \sqrt{\delta'_k}(\e'_k)^{n-1}\# \tilde{Z}'_k$ and consequently
\begin{equation}\label{controltZpk}
(\e'_k)^{n-1}\# \tilde{Z}'_k\ \le\ \sqrt{\delta'_k}\,.
\end{equation}
This controls the total surface of the squares ${C'_i}^k$ such
that in the corresponding cylinder $C_i^k$, the measure
of $E_k+\e'_kB(0,b)\setminus E_k$ is far from the measure of
a perfectly straight strip of width $\e'_k b$. In the other cylinders,
we will be able to show that the boundary of $E_k$ is almost flat.

We see at this point that \eqref{bigestimk} still holds if $Z_k$ is
replaced with $\tilde{Z}_k$, and $Z'_k$ with $Z'_k\cup \tilde{Z}'_k$. Together
with~\eqref{eq:BUblocal}
(with again $A=C(0,1-\beta)$) it follows that
\begin{equation*}
\limsup_{k\to \infty} \frac{1}{\e'_k} \int_{C(0,1-\beta)\setminus
\bigcup_{i\in\tilde{Z}_k} C_i^k} \left|\chi_{E_k+\e'_k B(0,b)}-\chi_{E_k}\right|\,dy
\ =\ 0\,,
\end{equation*}
and as a consequence
\begin{equation}\label{ErrorZero}
\limsup_{k\to \infty} \frac{1}{\e'_k} \int_{C(0,1-\beta)\setminus
\bigcup_{i\in\tilde{Z}_k} C_i^k} \left|\chi_{E_k+\e'_k C}-\chi_{E_k}\right|\,dy
\ =\ 0\,.
\end{equation}

We now need to estimate the quantity $(1/\e'_k)\int_{C_i^k}|\chi_{E_k+\e'_kC}
-\chi_{E_k}|\,dy$ for $i\in \tilde{Z}_k$, hence when
\begin{equation}\label{estimgood1}
b|{C'_i}^k|\left(1- K\sqrt{\delta_k}\right)
\ \le\
\frac{1}{\e'_k}\int_{C_i^k}\left|\chi_{E_k+\e'_k B(0,b)}-\chi_{E_k}\right|\,dy
\ \le\ b|{C'_i}^k|\left(1- K\sqrt{\delta_k}+
\frac{\sqrt{\delta'_k}}{b\theta^{n-1}}\right)
\end{equation}
and (from~\eqref{controlsetlwbd})
\begin{equation}\label{estimgood2}
|{C'_i}^k\setminus D_i^k|\le 2|{C'_i}^k|
\frac{\sqrt{\delta_k}}{\theta^{n-1}\sqrt{\beta}}\,.
\end{equation}
The estimate will rely on the fact that,
whenever~\eqref{estimgood1}-\eqref{estimgood2} hold,
the boundary of $E_k$ must be flat enough so that we can control
also the volume of $(E_k+\e'_kC) \setminus E_k$.

We choose $k\in\N$ and $i\in \tilde{Z}_k$ so that \eqref{estimgood1} and
 \eqref{estimgood2} hold, and consider the change of variable
$y= z^k_i+\e'_k \hat{y}$. We let $F=(E_k-z^k_i)/\e'_k$,
 $Q=(C_i^k-z_i^k)/\e'_k \supset (0,\theta)^{n-1}\times (-\sqrt{\beta}/\e'_k,
\sqrt{\beta}/\e'_k)$, $D=(D_i^k-z_i^k)/\e'_k$.
We find that
\begin{equation}\label{estimgood3}
\int_{Q} \left|\chi_{F+B(0,b)}-\chi_{F}\right|\,d\hat{y}
\ \le\ b \theta^{n-1}\left(1- K\sqrt{\delta_k}+
\frac{\sqrt{\delta'_k}}{b\theta^{n-1}}\right)\,,
\end{equation}
\begin{equation}\label{estimgood4}
\int_{(0,\theta)^{n-1}\times \left(-\frac{\sqrt{\beta}}{\e'_k},\frac{\sqrt{\beta}}{\e'_k}\right)}
\left|\chi_{F}-\chi_{\{\hat{y}_n\le 0\}}\right|\,d\hat{y}
\ \le\ a_i^k (\e'_k)^{-n}
\ \le\ \frac{\sqrt{\delta_k}}{\e'_k}\,,
\end{equation}
while \eqref{estimgood2} yields
\begin{equation}\label{estimgood22}
|(0,\theta)^{n-1}\setminus D|\le 2\sqrt{\frac{\delta_k}{\beta}}\,,
\end{equation}
where we recall that
\[
D\ =\ \left\{
\hat{y}'\in (0,\theta)^{n-1}\,:\,
\hs^1((\{\hat{y}'\}\times \R)\cap Q\cap F)>0\,,
\hs^1((\{\hat{y}'\}\times \R)\cap Q\setminus F)\ge b
\right\}\,.
\]

Let us set
$F_b = F+B(0,b) = \{\hat{y}\in\R^n\,:\, \dist(\hat{y},F)\le b\}$, and
$$
F_s = \{ \hat{y}\in F_b\,:\, \dist(\hat{y},\partial F_b)\ge b-s\}.
$$
We observe that $F\subset F_s$ for any $s\in [0,b]$.

By construction, given any point  $\hat{y}\in \partial F_s$, there
is a ball $B(\hat{z},b-s)$ such that $B(\hat{z},b-s)\cap F_s=\emptyset$
and $\hat{y}\in\partial B(\hat{z},b-s)$. Formally, it means that
the curvature of $\partial F_s$ is less than $1/(b-s)$. However, a similar
inner ball condition (with radius $s$) is not guaranteed.
We introduce the set $\tilde{F}_s$, which is the union of all balls
of radius $s$ which are contained in $F_s$. We have for $s\in (0,b)$
\begin{equation}
F+B(0,s)\, =\, \{\hat{y}\in\R^n\,:\, \dist(\hat{y},F)\le s\}
\ \subseteq\ \tilde{F_s}\ \subseteq\ F_s\,.
\end{equation}
We will show that if $k$ is large enough,  the boundary of $F_s$
is essentially flat inside $Q$.
Let us first establish that the boundary $\partial F_s$ crosses ``most''
of the vertical lines in the cylinder $(0,\theta)^{n-1}\times \R$.
We let:
\begin{align*}
D' &\ =\ \left\{ \hat{y}'\in (0,\theta)^{n-1}\,:\,
\hs^1((\{\hat{y}'\}\times \R_-)\cap Q \cap F) \ =\ 0\right\}\,,
\\
D'' &\ =\ \left\{ \hat{y}'\in (0,\theta)^{n-1}\,:\,
(\{\hat{y}'\}\times \R_+)\cap Q\ \subset F_b\right\}\,,
\\
D'''&\ =\ \left\{ \hat{y}'\in (0,\theta)^{n-1}\,:\,
\hs^1((\{\hat{y}'\}\times \R)\cap Q \cap (F_b\setminus F)) \ \ge\ 2b\right\}\,.
\end{align*}

The definition of $D$ ensures that if
$\hat{y}'\in D$,
$|(\{\hat{y}'\}\times \R)\cap Q \cap (F_b\setminus F)| \ \ge\ b$.
From~\eqref{estimgood3} and~\eqref{estimgood22}, we have that
\[
b\theta^{n-1}+\sqrt{\delta'_k}
\ \ge\ b |D\setminus D'''| + 2b |D'''|
\ \ge\ b |D| + b|D'''|
\ \ge\ b\theta^{n-1} + b|D'''| - 2b\sqrt{\frac{\delta_k}{\beta}}\,,
\]
hence
\begin{equation}\label{estimD3}
|D'''|\ \le\ \frac{1}{b}
\sqrt{\delta'_k}+2\sqrt{\frac{\delta_k}{\beta}}\,.
\end{equation}
We easily deduce from~\eqref{estimgood4} that both $|D'|$
and $|D''\setminus D'''|$ are bounded by a constant
times $\sqrt{\delta_k}$:
indeed, if $\hat{y}'\in D''\setminus D'''$,
$|(\hat{y}'\times \R_+)\cap Q\setminus F|\le 2b$ so that
\[
\left| \left(F\triangle \{ \hat{y}_n\le 0\}\right)
\cap
\left\{|\hat{y}_n|\le {\textstyle \frac{\sqrt{\beta}}{\e'_k} } \right\}
\right|
\ \ge\ \frac{\sqrt{\beta}}{\e'_k}-2b\,,
\]
and~\eqref{estimgood4} yields
\[
|D''\setminus D'''|
\ \le\ \frac{\sqrt{\delta_k}}{\e'_k}\frac{\e'_k}{\sqrt{\beta}-2b\e'_k}
\ \le\ 2\sqrt{\frac{\delta_k}{\beta}}
\]
as soon as $\e'_k\le \sqrt{\beta}/4b)$.
The estimate for $|D'|$ is even simpler.
It follows from this and~\eqref{estimD3}
that there exists a constant $K'$ (still depending on $\theta,\beta$)
such that
\begin{equation}\label{estimD1D2}
|D'\cup D''|\ \le\ |D'|+|D''\setminus D'''|+|D'''|\ \le\
K'\left(\sqrt{\delta_k}+\sqrt{\delta'_k}\right)\,.
\end{equation}

Now, each time $\hat{y}'\in (0,\theta)^{n-1}\setminus (D'\cup D'')$,
one can find
in $(\{\hat{y}'\}\times \R_-)\cap Q$ points which belong to $F$
(hence to $F_s$ or $\tilde{F_s}$ for all $s\in [0,b]$),
and in $(\{\hat{y}'\}\times \R_+)\cap Q$, points which are not in $F_b$:
as a consequence there are also, in  $(\{\hat{y}'\}\times \R)\cap Q$,
points in $\partial F_s$ or
$\partial\tilde F_s$, for any $s\in [0,b]$, with the set $F_s$ or $\tilde F_s$
``below'' the point (or the normal, if it exists, pointing upwards).
It follows that
\begin{equation}\label{estimlength}
\hs^{n-1}(\partial F_s\cap Q) \ \ge\ \theta^{n-1}-
K'\left(\sqrt{\delta_k}+\sqrt{\delta'_k}\right)\,,
\end{equation}
and (applying the coarea formula
to the distance function to $\partial F_b$) that
\begin{equation}\label{estimgood33}
|(F_b\setminus F_0)\cap Q|\ =\ \int_0^b \hs^{n-1}(\partial F_s\cap Q)\,ds
\ \ge\ b\theta^{n-1}
-bK'\left(\sqrt{\delta_k}+\sqrt{\delta'_k}\right)\,.
\end{equation}
Inequalities \eqref{estimgood3} and~\eqref{estimgood33} yield that
(here the constant $K'$ may vary from line to line, keeping the same
kind of dependency in the parameters)
\begin{equation}\label{estimrest}
|(F_0\setminus F)\cap Q| \ \le\
bK'\left(\sqrt{\delta_k}+\sqrt{\delta'_k}\right)\,.
\end{equation}

Let us choose $\eta>0$, small, and observe that (using~\eqref{estimlength}
and~\eqref{estimgood3})
\begin{multline}\label{estimpartiel}
\int_{\eta}^{2\eta} \hs^{n-1}(\partial F_s\cap Q)\,ds
\ =\ |(F_b\setminus F_0)\cap Q|-\int_{(0,b)\setminus (\eta,2\eta)}\hs^{n-1}(\partial F_s\cap Q)\,ds
\\ \le\
b \theta^{n-1}+\sqrt{\delta'_k}-(b-\eta)\left(
\theta^{n-1} - K'\left(\sqrt{\delta_k}+\sqrt{\delta'_k}\right)\right)
\ \le\ \eta\theta^{n-1} + K' \left(\sqrt{\delta_k}+\sqrt{\delta'_k}\right)
\end{multline}
so that there exists $\bar s\in [\eta,2\eta]$ with
\begin{equation}\label{flatness}
 \hs^{n-1}(\partial F_{\bar s}\cap Q)\ \le\ \theta^{n-1}+\frac{K'}{\eta}
\left(\sqrt{\delta_k}+\sqrt{\delta'_k}\right)\,.
\end{equation}

The surface $\partial F_{\bar s}$ is ``almost'' like a $C^{1,1}$ graph,
and converges to a flat surface, with convergence of its measure
as $k\to \infty$. If it were a graph, it would be easy to deduce
uniform convergence. Let us show that in this setting, it must
also be essentially flat. More precisely, we will establish it for
$\partial\tilde{F}_{\bar s}$.

Consider first
a point $\hat{y}\in \partial F_{\bar s}\cap\partial \tilde{F}_{\bar s}$,
that is, where $F_{\bar s}$ has both an outer ball condition with
radius $(b-\bar s)$ and an inner ball condition with radius $\bar s$.
In particular, there are at $\hat{y}$ two
tangent balls to $\partial F_{\bar s}$ of radius $\eta$ inside and
outside the set. The common normal to these balls is normal to
$\partial F_s$ (and $\partial \tilde{F}_s$) and we denote it
$\bar\nu$. Given $\zeta>0$, let us assume that $|\bar\nu \cdot\nu|\le
1-\zeta$.

\begin{figure}[h]
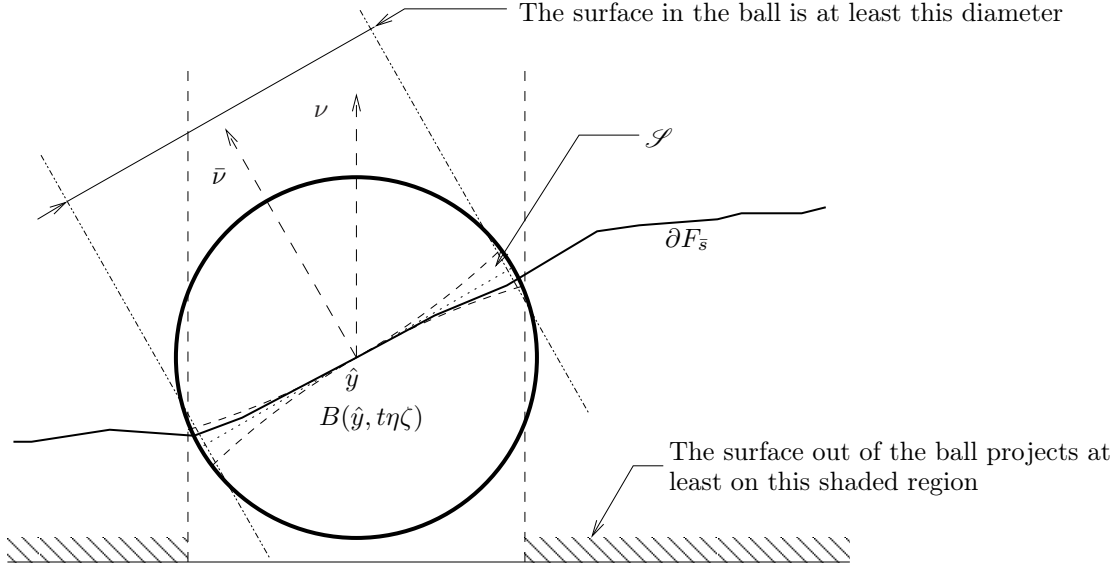

\begin{center}
\input ballFs.tex
\end{center}
\caption{If the normal at $\hat{y}$ to $\partial F_{\bar s}$ is away from $\nu$,
then its surface must exceed $\theta^{n-1}$ by
some quantity which is estimated.}
\label{FigBall}
\end{figure}

Then, for $t$ small, we consider the ball $B(\hat{y},t\eta\zeta)$ (which
we assume is in $Q$, and we let $r:=t\eta\zeta$).
The surface $\partial F_{\bar s}$ passes,
in $B(\hat{y},r)$ in between two spherical caps of radius $\eta$,
which are tangent in $\hat{y}$ and normal at that point to $\bar\nu$.
We call $\mathscr{S}$ the subset
of $B(\hat y,r)$ bounded by these two caps (see Figure~\ref{FigBall}).
A simple calculation
shows that the trace of these spherical caps on
the sphere $\partial B(\hat{y},r)$
is given by the intersection of this sphere with the
hyperplanes $\{(y-\hat{y})\cdot\bar \nu=\pm \frac{t\zeta}{2}r\}$
(hence,
$\mathscr{S}\subset B(\hat{y},r)
\cap \{|(y-\hat{y})\cdot\bar\nu|<\frac{t\zeta}{2}r\}$).
In particular, the surface $\hs^{n-1}(\partial F_{\bar s}\cap B(\hat{y},r))$
can be estimated from below with the surface of the corresponding
discs, that is, $\alpha_n r^{n-1}\sqrt{1-t^2\zeta^2/4}^{n-1}$.

Let us now estimate from below the surface
$\hs^{n-1}(\partial F_{\bar s}\cap Q\setminus B(\hat{y},r))$. Since we
know that given any $\hat{y}'\in (0,\theta)^{n-1}\setminus (D'\cup D'')$,
$\partial F_{\bar s}\cap (\{\hat y'\}\times \R)\cap Q\neq\emptyset$,
it is enough to estimate from above the projection of $\mathscr{S}$
onto $(0,\theta)^{n-1}$, which we denote by $\Pi_\nu (\mathscr{S})$.
This, in turn, is bounded by the projection of
\[
B(\hat{y},r) \cap
\left\{y\,:\,|(y-\hat{y})\cdot\bar\nu|<\frac{t\zeta}{2}r\right\}
\ =\
\left\{ \hat{y}+ r (s\bar\nu + \xi)\,:\, |s| < \frac{t\zeta}{2}\,,
|\xi|\le \sqrt{1-s^2}\,, \xi\cdot\bar\nu=0\right\}\,.
\]
Now, this projection is a subset of the vertical projection
of the diameter of $B(\hat{y},r)$ perpendicular to $\bar\nu$,
that is, $\Delta=\{\hat{y}+r\xi\,:\, |\xi|\le 1\,,\xi\cdot\bar\nu=0\}$,
plus the disk $\Pi_\nu(B(0,rt\zeta/2))$.
It follows (see the expansion of the volume of Minkowski sums of convex
sets in~\cite{Schneider}, {\it cf} Remark~\ref{remConv}) that
\[
|\Pi_\nu(\mathscr{S})|\ \le\ |\Pi_\nu(\Delta)|
\,+\,\P(\Pi_\nu(\Delta))\frac{rt\zeta}{2}
+o(rt\zeta) \,.
\]
Here, $\P(\Pi_\nu(\Delta))$ is the $(n-2)$-dimensional perimeter
of $\Pi_\nu(\Delta)$ in $(0,\theta)^{n-1}$, and
a simple scaling argument shows that $o(rt\zeta)$ is
of the form $r^{n-1}o(t\zeta)$, where the latter `$o$'
depends only on the geometry of the vertical projection of the unit
ball, that is, on $\nu\cdot\bar\nu$ --- and, in fact, would be largest for
$\nu=\bar\nu$.
Now, since $\P(\Pi_\nu(\Delta))\le\hs^{n-2}(\partial\Delta) = (n-1)\alpha_n r^{n-2}$, we obtain that
\[
|\Pi_\nu(\mathscr{S})|\ \le\ |\nu\cdot\bar\nu|\alpha_n r^{n-1}
\,+\, 2(n-1)\alpha_n r^{n-2} \frac{rt\zeta}{2}
\ \le\ \alpha_n r^{n-1}(1-(1-(n-1)t)\zeta)
\]
if $t$ is small enough.
It follows
\begin{equation}\label{estimbelowFs}
\hs^{n-1}(\partial F_{\bar s}\cap Q)
\ \ge\ \theta^{n-1}-K'(\sqrt{\delta_k}+\sqrt{\delta_k'})
+\alpha_n r^{n-1}
\left(\sqrt{1-\frac{t^2\zeta^2}{4}}^{n-1} -1 + (1-(n-1)t)\zeta \right)\,.
\end{equation}
For $t=0$, the quantity between the right-hand side parentheses is $\zeta>0$,
and it decreases with $t$. It follows that one can find $\bar t>0$ (depending
only on $n$ and $\zeta$)
such that \eqref{estimbelowFs} reads
 \[
 \hs^{n-1}(\partial F_{\bar s}\cap Q)\ \ge\ \theta^{n-1}
 \,-\,K'(\sqrt{\delta_k}+\sqrt{\delta'_k})
 \,+\, \alpha_n(\bar t\eta\zeta)^{n-1}\frac{\zeta}{2}\,.
 \]
Together with~\eqref{flatness},
it follows that if $k$ is large enough (depending on $K',\eta,\zeta$),
we get a contradiction, and therefore
$|\bar \nu\cdot\nu| > 1-\zeta$, provided $\hat{y}$
is at distance at least $\bar t\eta\zeta$ from $\partial Q$.

We must observe at this point that we also have $\bar\nu\cdot\nu\ge 0$.
Indeed, the same proof will show that if, for instance,
$\bar\nu\cdot\nu\le -1/2$, then for $k$ large enough, \eqref{flatness}
cannot hold. Indeed, in this case, thanks to~\eqref{estimgood4},
the surface of $\partial F_{\bar s}$ near $\hat{y}$, which is of
order $\eta^{n-1}$, has to be added to a surface of $\partial F_{\bar s}\cap Q$
(out of $(D'\cup D'')\times \R$)
already of order $\theta^{n-1}$, a contradiction if $k$ is large enough.
It follows than when $k$ is large, one must have
$\bar \nu\cdot\nu > 1-\zeta$, at any
$\hat{y}\in  \partial F_{\bar s}\cap \partial \tilde{F}_{\bar s}\cap Q$,
at distance at least $\bar t\eta\zeta$ from $\partial Q$.

We can deduce that $\partial \tilde{F}_{\bar s}$ is almost flat.
The reason is the following: given $B(\hat{y},\bar s)\subset \tilde{F}_{\bar s}$,
if we translate this ball in any direction
$(1-\zeta)e-\sqrt{\zeta(2-\zeta)}\nu$, for $e$ a unit horizontal vector (normal
to $\nu$), then it will never touch $\partial F_{\bar s}$, at least
until it reaches a distance $\bar t\eta\zeta$ to $\partial Q$. Otherwise,
necessarily, it would touch at some point where
$\bar \nu\cdot [(1-\zeta)e-\sqrt{\zeta(2-\zeta)}\nu]\ge 0$,
which yields $\bar\nu\cdot \nu\le 1-\zeta$. We denote
\[
Q^\zeta\ =\
\left\{\hat{y}\in Q\,:\,\dist(\hat{y},\partial Q)\ge \bar t\eta\zeta\right\}\,.
\]
Since by construction each point in $\partial \tilde{F}_{\bar s}$ belongs
to the boundary of a ball $B(\hat{y},\bar s)\subset \tilde{F}_{\bar s}$,
we find as a consequence
(taking also into account \eqref{estimgood3} and~\eqref{estimgood4})
that in $Q^\zeta$,
$\tilde{F}_{\bar s}$ is the subgraph $\{\hat{y}_n\le u(\hat{y}')\}$
of a Lipschitz function $u$, with Lipschitz constant $\sqrt{\zeta(2-\zeta)}/(1-\zeta)\le 2\sqrt{\zeta}$ (since $\zeta$ is small).

We deduce that there exists a value $\sigma$ such that
\begin{equation}\label{FFlat}
\partial \tilde F_{\bar s} \cap Q^\zeta\ \subset
 \{ \hat{y}\in Q\,:\,\sigma \ \le\ \hat{y}\cdot \nu \le\ \sigma+2\sqrt{\zeta}\theta \}
\end{equation}
for $k$ large enough, with moreover
$F\cap Q^\zeta\subset \tilde F_{\bar s}\cap Q^\zeta\subset
\{ \hat{y}\in Q^\zeta\,:\,\hat{y}\cdot \nu \le\ \sigma+2\sqrt{\zeta}\theta \}$.
In particular, it follows that $((F\cap Q^\zeta)+C)\cap Q^\zeta\subset
\{ \hat{y}\cdot \nu \le\ h_C(\nu)+\sigma+2\sqrt{\zeta}\theta \}$.
Observe however that this does not control the volume of the possible
points in $(F+C)\cap Q$ which could come from $(F\setminus Q^\zeta)+C$.
Since $C\subset B(0,b)$ and we can assume $\bar t\eta\zeta\le b$,
these points are outside of $(2b,\theta-2b)^{n-1}\times\R$, so that
\begin{multline}\label{finalsplit}
((F+C)\setminus F) \cap Q
\ \subseteq\  ((F_{\bar s}\setminus F)\cap Q)\,\cup\,
(((F+C)\setminus F_{\bar s})\cap Q)
\\ \subseteq ((F_{\bar s}\setminus F)\cap Q) \,\cup\,
((F_b\setminus F_{\bar s})\cap Q)\setminus ((2b,\theta-2b)^{n-1}\times\R))
\\
\hspace*{4.6cm}\cup ((((F\cap Q^\zeta)+C)\setminus F_{\bar s})\cap Q^\zeta
\\ \subseteq ((F_{\bar s}\setminus F)\cap Q) \,\cup\,
((F_b\setminus F_{\bar s})\cap Q)\setminus ((2b,\theta-2b)^{n-1}\times\R))
\\
\cup\,
\left\{
 \hat{y}\in Q^\zeta\,:\, \sigma\le \hat{y}\cdot \nu
\le\ \sigma+2\sqrt{\zeta}\theta+h_C(\nu) \right\}\,.
\end{multline}

The last set in the right-hand side has volume bounded
by $\theta^{n-1}(h_C(\nu)+2\sqrt{\zeta}\theta)$, which is the desired order, and we need to show that the two other
sets are much smaller.
To estimate the volume of the second set, we first check
that exactly for the same reasons for which~\eqref{estimlength} holds,
we have for $s\in [0,b]$
\begin{equation*}
\hs^{n-1}(\partial F_s\cap Q\cap ((2b,\theta-2b)^{n-1}\times \R)) \ \ge\ (\theta-4b)^{n-1}-
K'\left(\sqrt{\delta_k}+\sqrt{\delta'_k}\right)\,,
\end{equation*}
so that, still integrating from $0$ to $b$ and using the coarea formula,
\begin{equation*}
|(F_b\setminus F_0)\cap Q\cap ((2b,\theta-2b)^{n-1}\times \R)|
\ \ge\ b(\theta-4b)^{n-1}
-bK'\left(\sqrt{\delta_k}+\sqrt{\delta'_k}\right)\,.
\end{equation*}
Hence, using again~\eqref{estimgood3}, we find
\begin{multline}\label{estimborders}
|(F_b\setminus F_0)\cap Q\setminus ((2b,\theta-2b)^{n-1}\times \R)|
\\ \le\ b\theta^{n-1}+\sqrt{\delta'_k}-b(\theta-4b)^{n-1}
+bK'\left(\sqrt{\delta_k}+\sqrt{\delta'_k}\right)\,.
\\ \le\ 4(n-1)b^2\theta^{n-2} +
K''\left(\sqrt{\delta_k}+\sqrt{\delta'_k}\right)\,.
\end{multline}

Exactly in the same way as~\eqref{estimpartiel} we also see that
\[
|(F_{\bar s}\setminus F_0)\cap Q|\ \le\ \bar{s}\theta^{n-1}
\,+\,K'\left(\sqrt{\delta_k}+\sqrt{\delta'_k}\right)\,,
\]
which combined with~\eqref{estimrest} yields
\[
|(F_{\bar s}\setminus F)\cap Q|\ \le\ 2\eta \theta^{n-1}\,+\,
(b+1)K'\left(\sqrt{\delta_k}+\sqrt{\delta'_k}\right)\,.
\]
This and~\eqref{estimborders} show that~\eqref{finalsplit} can
be estimated as follows:
\[
|((F+C)\setminus F)\cap Q|\ \le\ \theta^{n-1}(2\sqrt{\zeta}\theta+h_C(\nu))
\,+\,  2\eta \theta^{n-1}\,+\, 4(n-1) b^2\theta^{n-2} + \mathscr{R}_k
\]
where $\mathscr{R}_k$ is a rest which goes to zero
with $\delta_k$ and $\delta'_k$, and does not depend on the particular
cylinder $C_i^k$ we were examining.
Returning to the original sets $C_i^k$,
we find that if $k,i\in \tilde{Z}_k$ and $k$ is large enough,
\begin{equation}\label{estimcik}
\frac{1}{\e'_k}\int_{C_i^k} |\chi_{E_k+\e'_kC}-\chi_{E_k}|\,dy
\ \le\ |{C'_i}^k|\left(h_C(\nu) + 2\sqrt{\zeta}\theta
+ 2\eta + \frac{4(n-1) b^2}{\theta} + \frac{1}{\theta^{n-1}}\mathscr{R}_k
\right)\,.
\end{equation}
Together with~\eqref{ErrorZero}, \eqref{estimcik} yields that
\begin{equation*}
\limsup_{k\to 0} \frac{1}{\e'_k} \int_{C(0,1-\beta)}
 \left|\chi_{E_k+\e'_k C}-\chi_{E_k}\right|\,dy
\ \le\ \alpha_n(1-\beta)^{n-1} \left(h_C(\nu) + 2\sqrt{\zeta}\theta
+ 2\eta + \frac{4(n-1) b^2}{\theta}\right)\,.
\end{equation*}
Sending first $\zeta$, then $\eta$ to zero and eventually $\theta$
to $+\infty$, and
using~\eqref{eq:BUC} and~\eqref{eq:BUblocal}, we deduce
\begin{equation*}
\alpha_n g(\bar x)\ =\ \lim_{k\to 0} \frac{1}{\e'_k} \int_{B(0,1)}
 \left|\chi_{E_k+\e'_k C}-\chi_{E_k}\right|\,dy
\ \le\ b\alpha_n(1-(1-\beta)^{n-1})
\,+\,\alpha_n(1-\beta)^{n-1}h_C(\nu)\,,
\end{equation*}
and letting then $\beta\to 0$ yields the desired inequality.

\end{document}